# Determining the open pit to underground transition: A new method

*D Whittle[1], M Brazil[2], P A Grossman[3], J H Rubinstein[4], D A Thomas[5],*

1.
MAusIMM, Graduate Researcher, Department of Mechanical Engineering, The University of Melbourne, Parkville, Vic 3010. Email: dwhittle1@student.unimelb.edu.au

2.
Associate Professor and Reader, Department of Electrical and Electronic Engineering, The University of Melbourne, Parkville Vic 3010. Email: brazil@unimelb.edu.au

3.
Senior Research Fellow, Department of Mechanical Engineering, The University of Melbourne, Parkville, Vic 3010. Email: peterag@unimelb.edu.au

4.
MAusIMM, Professor, Department of Mathematics and Statistics, The University of Melbourne, Parkville, Vic 3010. Email: rubin@ms.unimelb.edu.au

5.
MAusIMM, Professor, Head of Department and Associate Dean Research, Department of Mechanical Engineering, The University of Melbourne, Parkville, Vic 3010. Email: doreen.thomas@unimelb.edu.au

# ABSTRACT

Many Ore Reserves are harvested by a combination of open pit and underground mining methods. In these cases there is often material that could be mined by either method, and a choice has to be made. The area containing this material is referred to as the transition zone. Deciding where to finish the open pit and start the underground is referred to as the transition problem and it has received some attention in the literature since the 1980s.

In this paper we provide a review of existing approaches to the transition problem encompassing: graph-theory based optimisation employing an opportunity cost approach; heuristics and integer programming. We also present a novel opportunity cost approach, allowing it to take into account a crown pillar, and show how the new approach can be best applied through the unconventional application of an existing mine optimisation tool.

# INTRODUCTION

The widespread application of mathematical optimisation of pit outlines had its beginnings with the release of optimisation packages by Whittle Programming (now owned by Dassault Systemes Geovia) in the 1980s. The software includes several optimisation engines, including one based on the graph algorithm developed by Lerchs and Grossmann (1965). This algorithm, commonly referred to as the LG Algorithm, provides an exact method to optimise the outline of an open pit mine. Picard (1976) showed that the pit outline optimisation problem was equivalent to a maximum flow problem, meaning that a large number of existing maximum flow algorithms could be applied. Giannini et al. (1991) developed a software package using a maximum flow algorithm and others (for example Zhao and Kim (1992)) have proposed or devised methods that produce equivalent results to the LG Algorithm with the potential to do so in less time. However the LG Algorithm is still thought to be the most commonly used pit optimiser.

The LG Algorithm only optimises the pit outline, but in combination with other optimisation tools, a wide range of analysis and planning decisions can be supported (e.g.  Ramazan (2007), Stone et al. (2004), Whittle (2011) and Whittle (2014)). Bienstock and Zuckerberg (2010) have published a method to simultaneously optimise the pit outline and the schedule subject to a range of side constraints.

The optimisation of underground mine plans has received less attention, though heuristics for stope envelope optimisation have been developed (e.g. Alford (1995), Alford and Hall (2009), Thomas and Earl (1999), Bai et al. (2014) and Sandanayake et al. (2015)) and some are now available in commercial software. In addition, algorithms and software have been developed to optimise access to stopes from the surface to either minimise cost (e.g. Brazil and Grossman (2008), Brazil et al. (2013)) or maximise Net Present Value (NPV) (Sirinanda et al. (2015a), Sirinanda et al. (2015b)).

Various authors have addressed the issue of the combined optimisation of open pit and underground mines and some have focused in particular on the optimisation of the transition from open pit mining to underground mining, and that leads us to this paper, the rest of which is organised as follows:

- Current Approaches - summarises the various ways in which academics and practitioners have viewed the transition problem and the state of the art as to tools applied to optimisation and decision-making.
- New Opportunity Cost Approach – a new mathematical model that adapts an existing pit optimisation algorithm to include the physical and economic impacts of a crown pillar separating an open pit mine from an underground mine.
- Demonstration of the New Opportunity Cost Approach – implements the new method in the Whittle software.
- Practical Considerations – discusses how the new method can be incorporated into a wider transition planning workflow.

# CURRENT APPROACHES

Nilsson (1981) gave a good account of the factors that need to be considered when planning a transition from open pit to underground mining, including geotechnical interactions, economic considerations and operational challenges. He recognised that the optimal pit design changes if deeper parts of the orebody can be mined by an underground mine and his illustration of the difference is included in Figure 1.

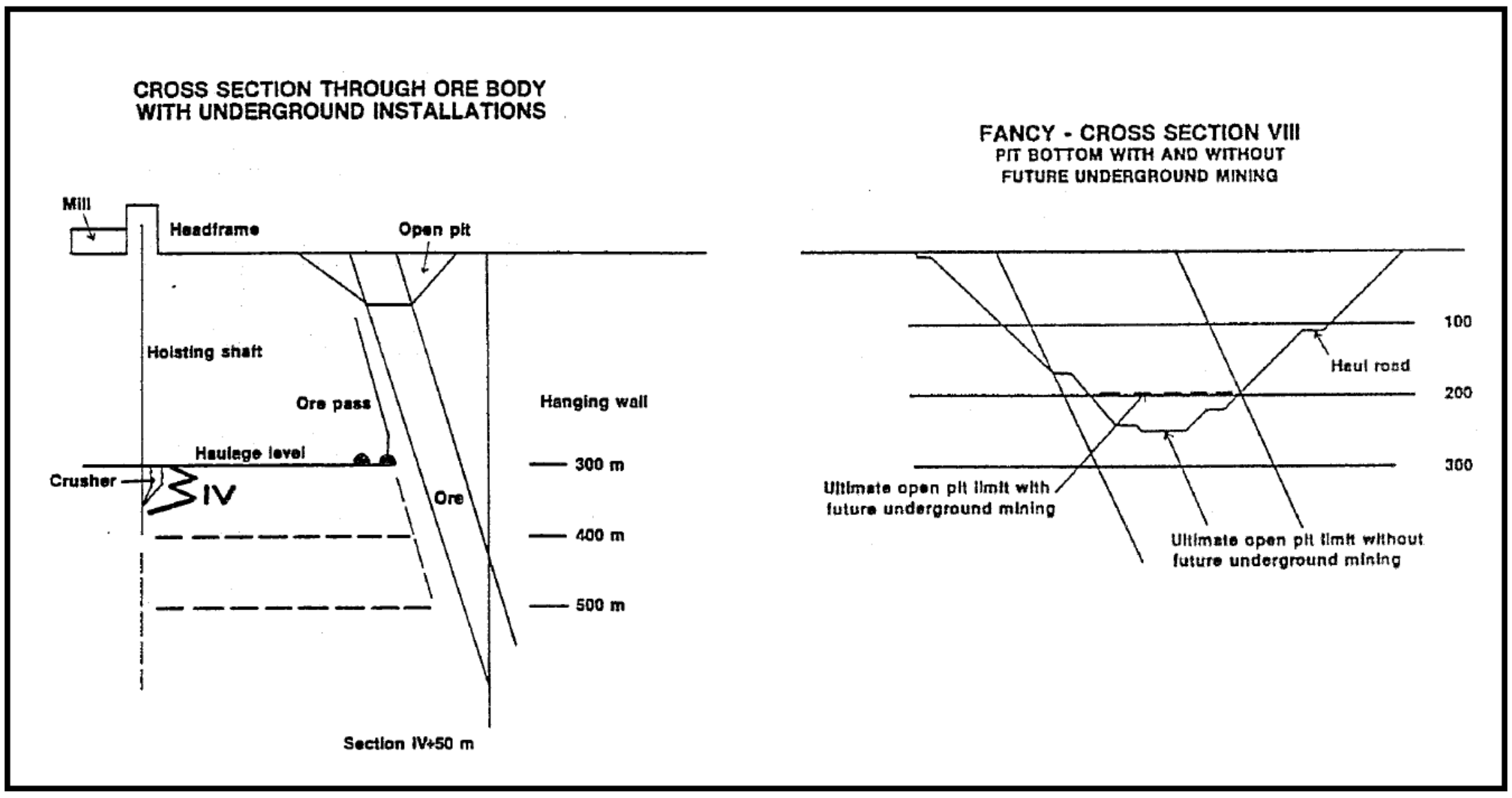

*Figure 1: Original illustrations from Nilsson (1981) showing the relationship between an open pit mine and an underground mine on the same orebody (left), and the typical change in the shape of the open pit mine as a result of the underground opportunity (right).*

Nilsson's paper predated the wide availability of pit optimisation software and he did not give a detailed account of the calculations leading to a change in the pit shape. Whittle (1990) incorporated a method into pit optimisation software that takes into account the value that ore has if mined by underground method. Consider a case in which some blocks can be mined by either open pit or underground method. For any such block, the value used for pit optimisation must be the difference between its open pit value and its underground value. The assumption underpinning this is that for a block that can be mined by either method, if it is not mined by open pit method, it will be mined by underground method. We will henceforth refer to this as the opportunity cost approach. Camus (1992) independently described an approach that will generate equivalent results.

When the opportunity cost approach is applied, the open pit almost always gets smaller. There is more than one way to generate a smaller pit using pit optimisation software, for example, it is common to use a technique called pit parameterisation to generate a family of pits by flexing the commodity price. However, the pit created using the opportunity cost approach may not match the shape and tonnage of any of the pits created using regular parameterisation techniques, due to the different ways in which the block values are calculated.

Chen et al. (2003) described a method along similar lines to Whittle (1990) and Camus (1992) but without exact optimisation, and stated that it is accepted practice in Russia and China. They included consideration of a crown pillar, which the earlier authors had not. A crown pillar is a body of rock left in place above the shallowest part of an underground mine to ensure stability in the ground above. The need for stability is driven by the land use, which in some cases is open cut mining. The crown pillar also acts to reduce or avoid the ingress of water to the underground mine and to ensure the stability of the cavity below. An indication of the position of a crown pillar is provided in Figure 2. When crown pillars are used, their design must take into account the geotechnical characteristics of the native rock and the planned sizes and shapes of the underground stopes (Carter (2000)). Note that the use of crown pillars is not universal, since caving of the ground above the underground mine is sometimes a desired or acceptable outcome, if good engineering and operational controls can be established.

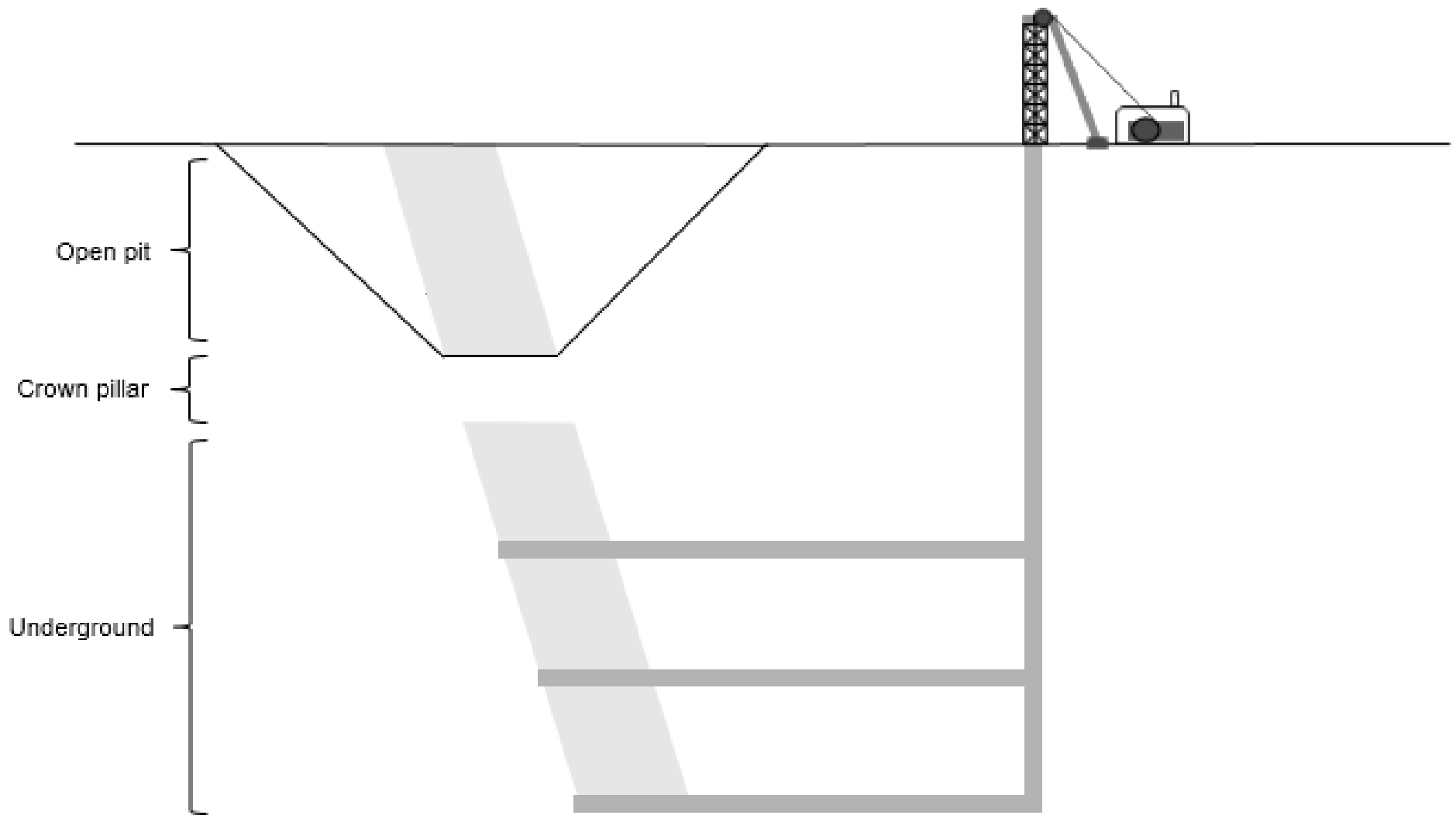

*Figure 2: Illustration of a crown pillar between an open pit and underground mine.*

What the above-mentioned methods have in common is the use of a pit optimisation tool (or a heuristic equivalent), which considers the cost associated with the lost opportunity to mine ore in the underground mine. The methods optimise only on the shape and size of the pit, and do so in order to maximise undiscounted cash flows. Pit optimisation can however be incorporated into a wider workflow with other optimisation tools to support a wide range of analysis and planning decisions: Finch and Elkington (2011) advocate for automated scenario analysis in which a number of candidate transition depths are evaluated with schedule optimisation software. Roberts et al. (2013) provide a case study using a number of different in-house and commercially available optimisation tools. The unnamed mine has an existing underground mine and is contemplating an open pit expansion that could mine material that would otherwise be mined in a future underground expansion. These authors were able to compare a number of different transition scenarios, each with optimised schedules and cut-offs.

Chung et al. (2015) formulated an integer programming model of the underground and open pit mines with crown pillar and made optimisation tractable by using very simple pit slope and stope envelope constraints. A number of researchers have also been working on optimisation approaches to take into account a crown pillar, and to optimise directly for the maximisation of Net Present Value (NPV), through the development of purpose built optimisation tools. Bakhtavar et al. (2008) published the first of their several papers on the topic, applying the opportunity cost approach through a heuristic algorithm. They used calculated stripping ratios (the ratio between waste mass and ore mass, with the former needing to be excavated in order to access the latter) rather than any optimisation for the determination of the open pit. Bakhtavar et al. (2009) determined the transition depth through the comparison of NPV calculations of both open pit and underground options for each level. Bakhtavar et al. (2012) applied integer programming in the consideration of transition depth, taking into account the crown pillar (with the tonnes removed from the underground). They compared a range of transition depths in order to determine the depth that gave the highest NPV, though applied only to a two dimensional model. Newman et al. (2013) formulated the transition as a large monolithic longest-path network problem. MacNeil and Dimitrakopoulos (2014) formulated the open pit to underground transition as a stochastic optimisation problem. In order to make the problem tractable, they represented the mines in a set of stratums: conceptually horizontal slices of the model. The authors incorporate consideration of a crown pillar, processing cut-off decisions and mining rate decisions in the optimisation model. The result is an assignment of stratums to open pit, crown pillar and underground mine respectively.

In summary, there are three general approaches:

- Use the LG Algorithm (or an equivalent) and an opportunity cost model to calculate the optimal transition point. The approach optimises net cash values, and cannot take account of a crown pillar.

- Use purpose-built optimisation tools employing mixed integer programming models that maximise NPV. In order to make the problem tractable, simplify the potential open pit and underground mining shapes considerably.
- Use workflows that rely on a combination of optimisation tools and judgement to make a wide range of planning decisions, including transition depth, in order to maximise NPV.

None of the approaches can handle the full range of complexities confronting mine planners. The LG Algorithm methods use a very capable and mature pit optimisation approach, but cannot model the crown pillar, or optimise for NPV directly. The purpose-built optimisation tools can handle NPV and the crown pillar, but are very simplistic in their treatment of the open pit and underground designs. In the next section we overcome some of these deficiencies, focusing especially on a new model for the application of the opportunity cost approach using the LG Algorithm that provides for very flexible modelling of the physical and economic characteristics of the crown pillar.

# NEW OPPORTUNITY COST APPROACH

We have devised a new opportunity cost approach through three modifications to the construction of the model presented to the LG Algorithm and a change to the interpretation of the result of the optimisation. We commence by restating the model developed by Lerchs and Grossmann (1965) and the conventional approach to applying the opportunity cost approach. We then describe the modifications.

## Lerchs Grossmann Model

The model Lerchs and Grossmann described was quite general. This restatement focuses on the mining industry implementation of the model (for example the use of regular rectangular blocks). We also update the mathematical notation. There are two related models:

1. Block model: A block model is comprised of regular rectangular blocks representing air, waste rock and mineralised material in a three-dimensional framework. Given price and cost information, it is possible to assign dollar values to each of the blocks.
2. Digraph: A digraph is comprised of vertices and arcs (directed connections between vertices). It is a defined type of model in a branch of mathematics called Graph Theory.

The open pit optimisation problem involves finding the subset of blocks in the block model that maximises total value whilst obeying pit slope constraints. However, in devising an optimisation algorithm, it is more convenient to implement it in a digraph, so we need a digraph representation of the block model.

Let the set $X$ of vertices $x_i$ and the set $A$ of arcs $a_k = (x_i, x_j)$ define a digraph $G = (X, A)$. In this digraph:

- Each vertex $x_i \in X$ corresponds to block $i$ in the block model, and has a weight $m_i \in \mathbb{R}$. The weight represents the value of the corresponding block $i$ in the block model. The weight $m_i = v_i - c_i$ where $v_i$ is the revenue and $c_i$ is the extraction cost. In the interests of brevity, we sometimes use the terms *vertex* and *block* interchangeably below.
- Each arc $a_k = (x_i, x_j)$, $a_k \in A$ represents a mining dependency, specifically relating to the need to uncover a block and to maintain maximum safe pit slopes. If mining the block represented by the vertex $x_i$, is dependent on the block represented by $x_j$ being mined, then the arc $a_k = (x_i, x_j)$ will be included in $A$.

A closed subgraph $G_Y$ of $G$ is a graph $(Y, A_Y)$, such that:

- $Y \subseteq X$, and if $x_i \in Y$ and $(x_i, x_j) \in A$, then $x_j \in Y$
- $A_Y \subseteq A$, and if $x_i \in Y$ and $x_j \in Y$ and $(x_i, x_j) \in A$, then $(x_i, x_j) \in A_Y$

Figure 3 includes an illustration of a two-dimensional model with arcs (left) and an example of a closed subgraph (right). In a two-dimensional model with square blocks, only three arcs per block are required to model 45-degree slopes. Thirty or more arcs per block are needed to model slopes accurately in a typical three-dimensional model.

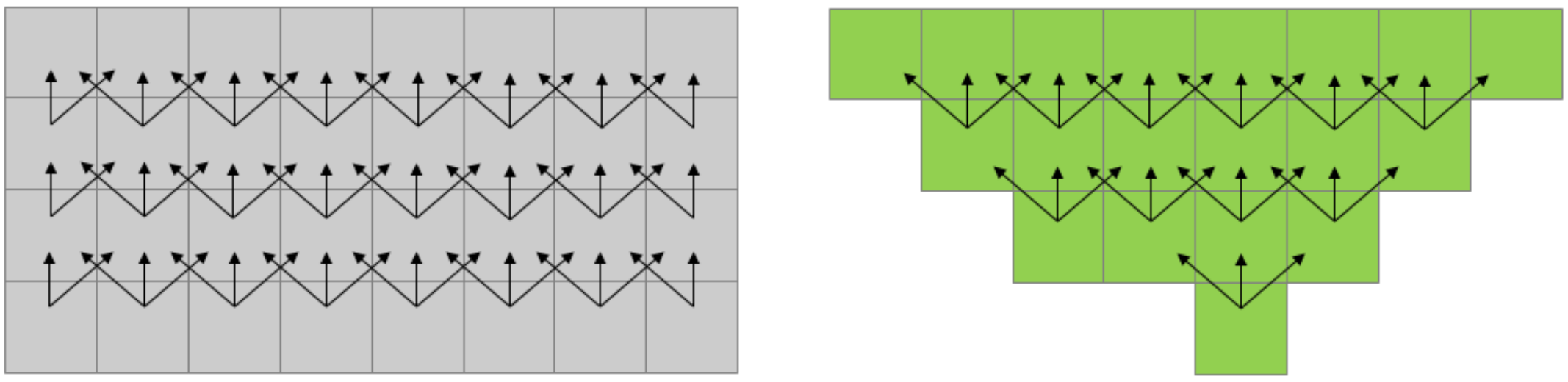

*Figure 3: A two-dimensional model with arcs representing block dependencies (left) and an example closed subgraph of the digraph (right).*

The abovementioned pit optimisation problem can now be restated in Graph Theory terms as finding a closed subgraph $G_Y$ of $G$ that maximises $M_Y$:

$$Max\ M_Y = \sum_{i:\ x_i \in Y} m_i$$

Lerchs and Grossmann (1965) defined a graph algorithm to solve the abovementioned problem. The algorithm has been widely accepted and applied and it is not necessary in this paper to restate it.

Figure 4 illustrates an optimal pit design. In this simple two-dimensional example, pit slopes are 45 degrees and there are just a few blocks. In practice, models are three-dimensional; have tens of thousands to millions of blocks; and have complex pit slope requirements. Lerchs and Grossmann's model and algorithm accommodate these practical requirements.

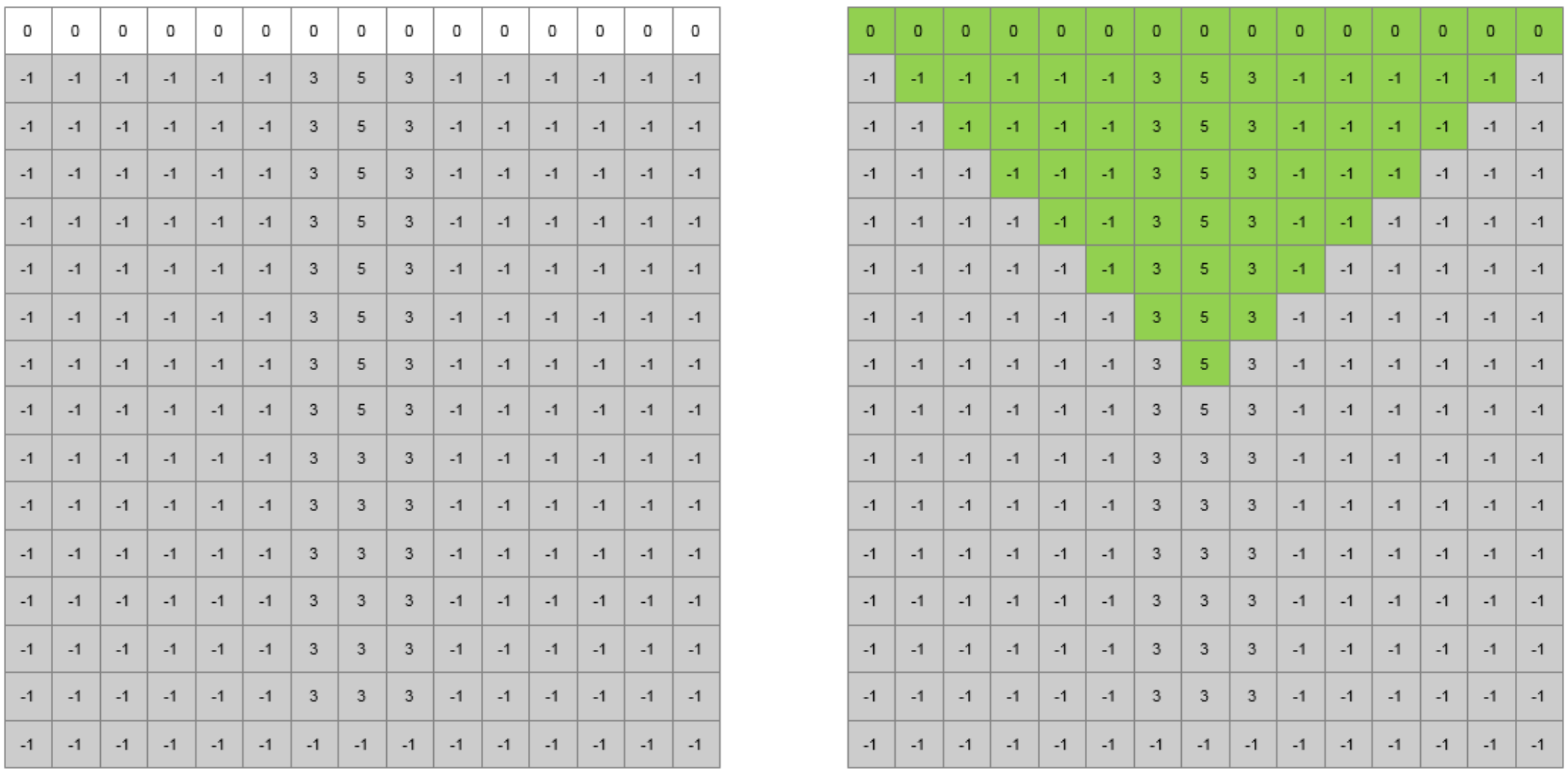

*Figure 4: Left: A simple model for open pit optimisation. Numbers represent values ($m_i$). Right: The optimal pit.*

## Conventional Opportunity Cost Approach

The conventional way to apply the opportunity cost approach (Whittle (1990), Camus (1992)) is to modify the way $m_i$ is calculated for each block:

$$m_i = v_i^p - c_i^p - (v_i^u - c_i^u)$$

Where:

$v_i^p$ The open pit revenue

$c_i^p$ The open pit extraction cost

$v_i^u$ The underground mining revenue

$c_i^u$ The underground mining extraction cost

There are differences in the way in which the open pit extraction cost $c_i^p$ and its underground counterpart $c_i^u$ are calculated. The open pit extraction cost is calculated assuming that the block is already uncovered and ready to mine. This assumption is important because the total value $M_Y$ of the closed subgraph $G_Y$ (representing the pit) includes the values for all the blocks that must be mined in accordance to the precedences defined by the arcs. Accordingly, if, in calculating $c_i^p$ for a block, the cost of uncovering that block is counted, then this leads to double counting of the uncovering cost. In contrast, the precedences for the blocks that will be mined by underground method (if not mined by open pit method) are not represented in the model. Accordingly, in order to calculate the underground extraction cost $c_i^u$ for a block it is necessary to include the costs to gain access, especially the prorated cost of excavating drifts, shafts and/or declines that are not explicitly represented in the model. The manner in which $c_i^u$ is calculated has important implications as to the practical application of the opportunity cost approach, discussed later in this paper (See Practical Considerations section below).

Figure 5 provides an illustration of how the opportunity cost approach works. The model as for Figure 4 is shown on the left. Underground values are shown in the centre, and the combined model and optimal pit with opportunity cost applied is shown on the right.

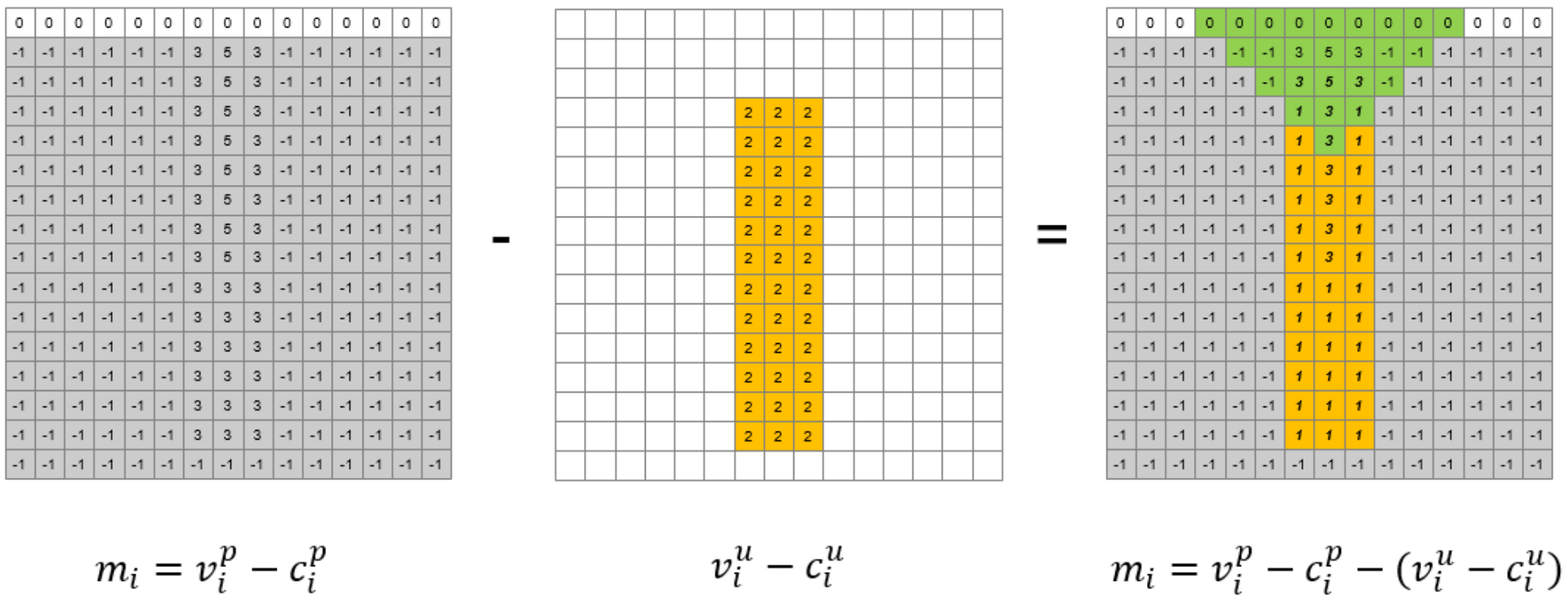

*Figure 5: Illustration of how the opportunity cost approach works.*

## First Modification – Alternative Opportunity Cost Approach

This first modification allows for the implementation of the opportunity cost approach in a different way. This new way produces identical results to the original, and provides a foundation for further improvements.

We introduce a second block model with the same dimensions as the first, but exclusively representing the blocks in potential underground stopes. Both the new underground model and the original open pit model occupy the same three-dimensional space.

To accommodate representation of this second model in the digraph we rename the original model's vertices from $x_i \in X$ to $x_i^p \in X^p$ and define new vertices $x_i^u \in X^u$ corresponding to the new underground mine block model.

$X = X^p \cup X^u$

Weights for vertices $x_i^p \in X^p$ are given by $m_i^p = v_i^p - c_i^p$. Weights for vertices $x_i^u \in X^u$ are given by $m_i^u = -(v_i^u - c_i^u)$.

We allow for new types of arcs by renaming the original set of arcs $A$ to $B$ mapping $X^p$ into $X^p$. Each arc $a_i^B \in B$ represents a mining dependency to maintain maximum safe pit slopes as illustrated previously in Figure 3.

In order to reproduce the conventional opportunity cost approach in the new digraph we introduce a set of arcs $C$ mapping $X^p$ into $X^u$. Each arc $a_i^C \in C$ connects vertex $x_i^p$ to its corresponding vertex $x_i^u$, for every $x_i^u$ that has a defined weight. The position of $x_i^u$ in three dimensional space is identical to that of $x_i^p$.

The effect is that if vertex $x_i^p \in X^p$ is included in $Y$, and if $a_i^C = (x_i^p, x_i^u)$ exists, then $x_i^u \in X^u$ is also included. This brings the opportunity cost into the calculation $M_Y = \sum_{i:\, x_i^p \in Y} m_i^p + \sum_{j:\, x_j^u \in Y} m_j^u$.

When the LG Algorithm is applied to this model, it finds a subgraph $G(Y) = (Y, A_Y)$ of $G = (X, A)$ that maximises $M_Y = \sum_{i:\, x_i^p \in Y} m_i^p + \sum_{j:\, x_j^u \in Y} m_j^u$ where $X = X^p \cup X^u$, and $A = B \cup C$. The new model reduces to the same type of graph problem as before, but the results are interpreted in a new way as illustrated in Figure 6. The open pit model is on the left and the underground model is on the right. For illustrative clarity, only 12 of the 36 arcs $a_i^C \in C$ are shown. Recall these arcs connect vertices $x_i^p$ to vertices $x_i^u$, for every vertex $x_i^u$ that has a defined weight. Green blocks on the left ($X^p \cap Y$) are blocks in the optimal open pit. Green blocks on the right ($X^u \cap Y$) are blocks no longer available for underground mining.

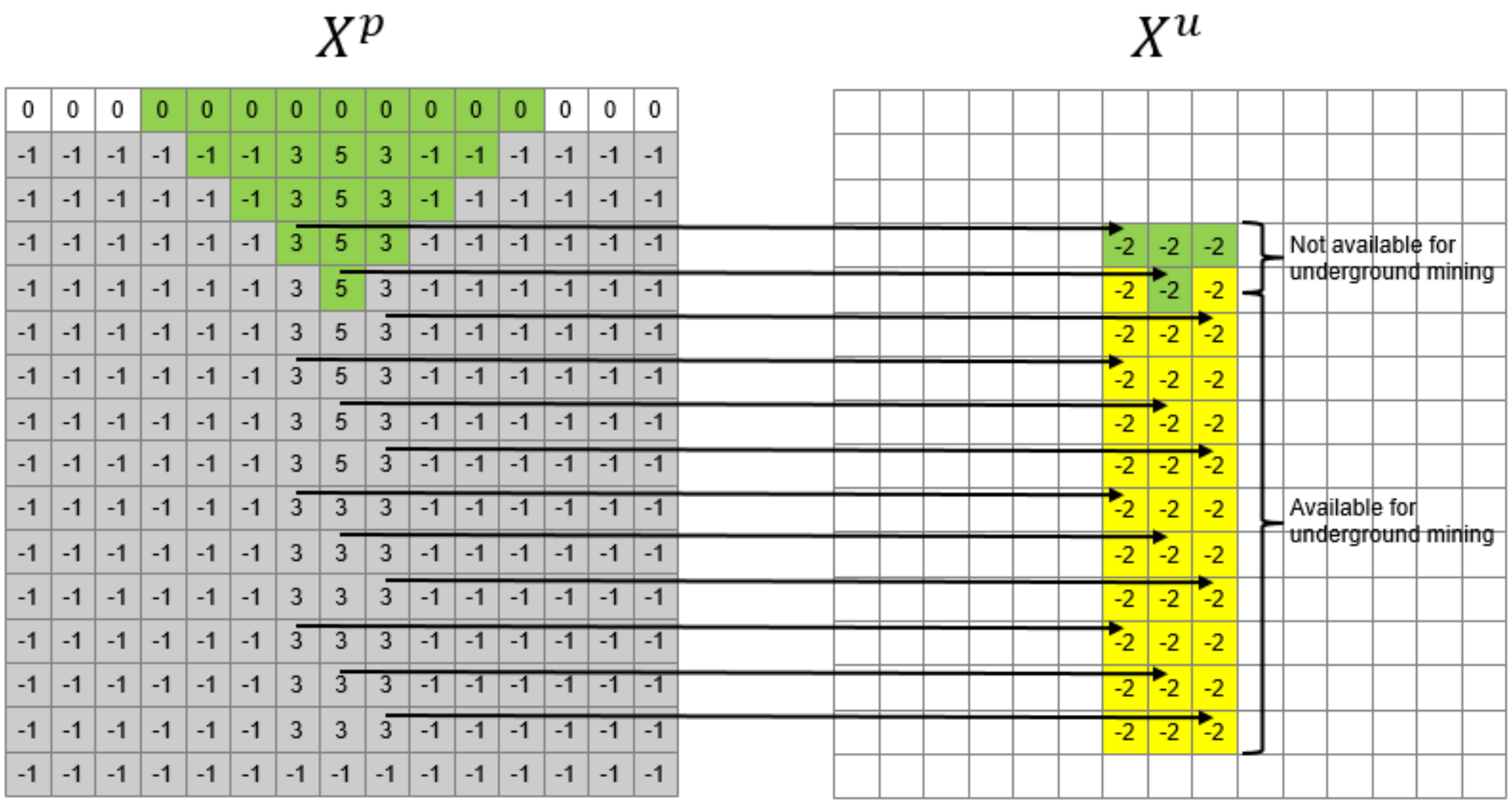

*Figure 6: Open pit model (left) and underground model (right).*

## Second modification – Simple crown pillar

A "simple" crown pillar is one in which the contour of the base of the pit is projected onto the top of the underground mining envelope with the projection vertically displaced by the required thickness of the crown pillar. "Simple" in this case refers to the formulation of the problem, rather than the shape of the mining envelopes.

The modification for a simple crown pillar is to change the way arcs $C$ map from $X^p$ into $X^u$. Each arc $a_j^C \in$ C connects vertex $x_i^p$ to vertex $x_j^u$, for every $x_j^u$ that has a weight defined. $x_j^u$ is a block under $x_i^u$, where the vertical separation is equal to the required thickness of the crown pillar and is given according to geotechnical constraints.

Figure 7 illustrates how this modification changes the outcome. In this example, the required height of the crown pillar is equal to twice the height of a block. Accordingly, arcs map from the open pit model $X^p$ into blocks two levels lower in the underground model $X^u$ (Compare with Figure 6). The additional blocks required for the crown pillar are shown. For illustrative clarity, only 12 of the 36 arcs $a_i^C \in C$ are shown.The blocks in the crown pillar can be identified as the set $P = \{x_i^u \in X^u \cap Y : x_i^p \notin Y\}$.

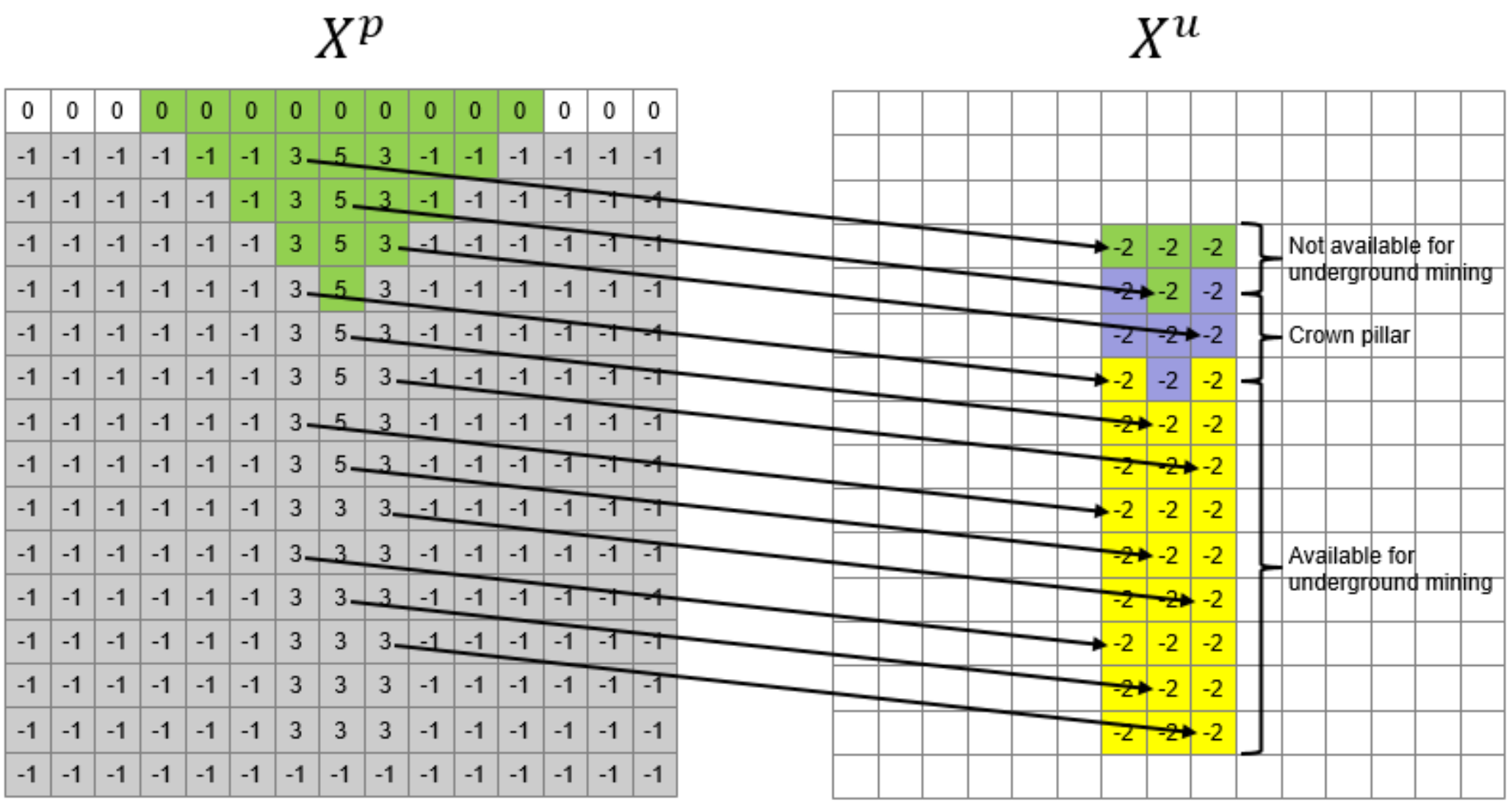

*Figure 7: Offset arcs take account of the crown pillar.*

In the example shown in Figure 6, the optimal pit size $(X^p \cap Y)$ has not changed, but other cases can be constructed in which the optimal pit size decreases. It is easy to see that this can be the case since the additional blocks required for the crown pillar reduce the value of the optimal pit. However, there are also cases in which the addition of a crown pillar increases the size of the optimal pit, and the Demonstration section below provides an example.

This "simple crown pillar" approach has achieved vertical separation between the open pit and the underground mine, but in the example shown in Figure 7, the shape of the crown pillar is almost certainly not practical for the underground mine. The third modification described below addresses this issue.

## Third modification – Well-formed crown pillar

In order to provide for a more practical shape to the crown pillar, we define a third set of arcs $D$ mapping $X^u$ into $X^u$. Arcs $a_i^{\mathrm{D}} \in D$ form directed cycles at each level comprised of arcs $\left(x_i^u, x_j^u\right)$, $x_i^u \neq x_j^u$ . Each cycle represents the required shape for the crown at a given level. If $x_i^p$ is included in $Y$ and if an arc $\left(x_i^p, x_j^u\right)$ exists, then $x_j^u$ will also be included in $Y$, along with every other vertex in its cycle. Figure 8 shows two examples. Each cycle defines the required shape for the top of the underground mining block. If any block is included in a closed subgraph $Y$, then every block in its cycle will also be included in $Y$. The example on the left would provide a flat top. The example on the right would provide an approximate dome-shaped top at higher levels, and flat tops at lower levels. These illustrations are two-dimensional for clarity. In a three-dimensional model, cycles can be used to create any non-overlapping three-dimensional shapes.

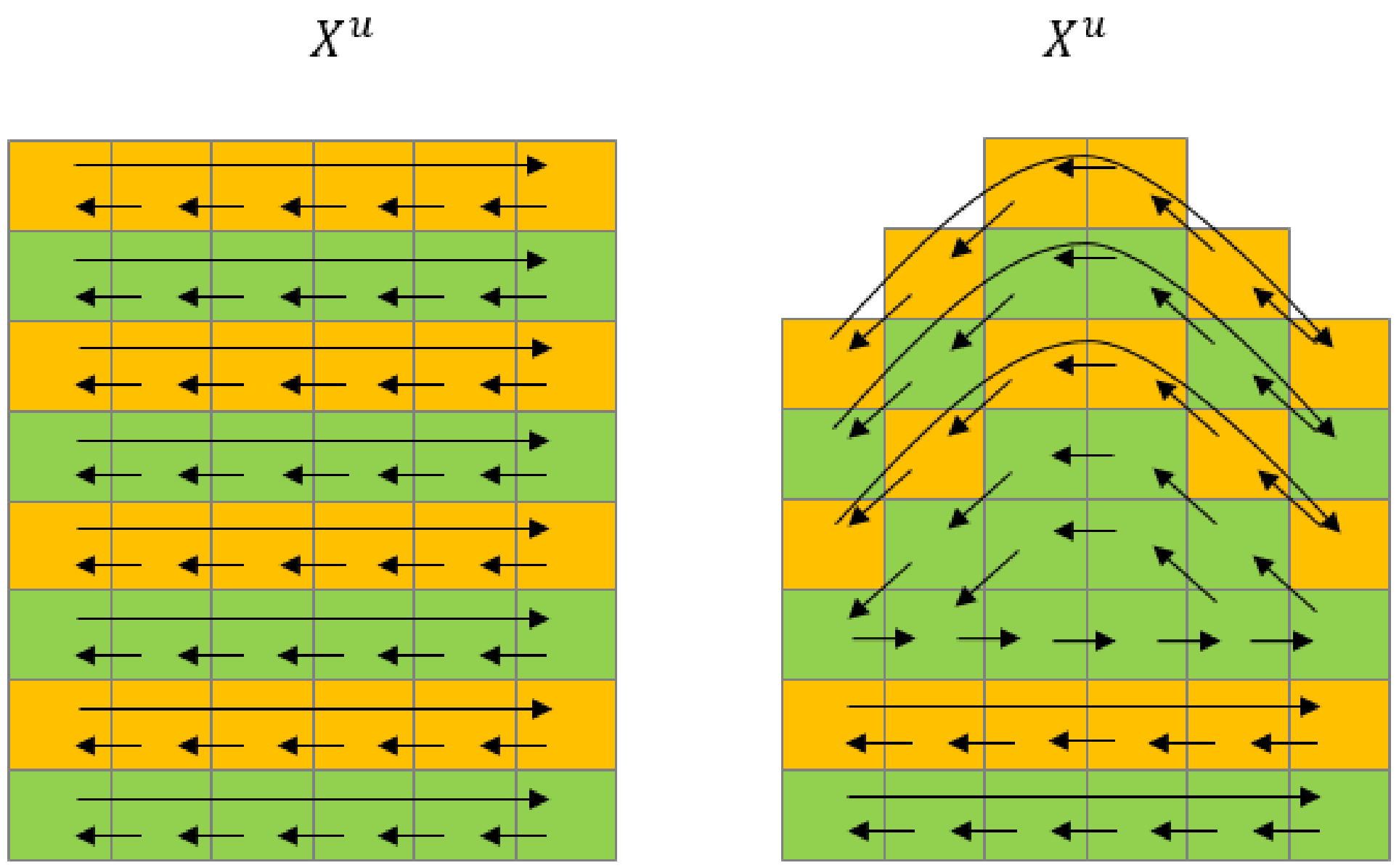

*Figure 8: Each cycle defines the required shape for the top of the underground mining block.*

Figure 9 provides an illustration. Arcs $a_i^D \in D$ (not shown) impose a flat top to the crown (see Figure 8 Left). In this case, the bottom of the pit has been flattened as well.

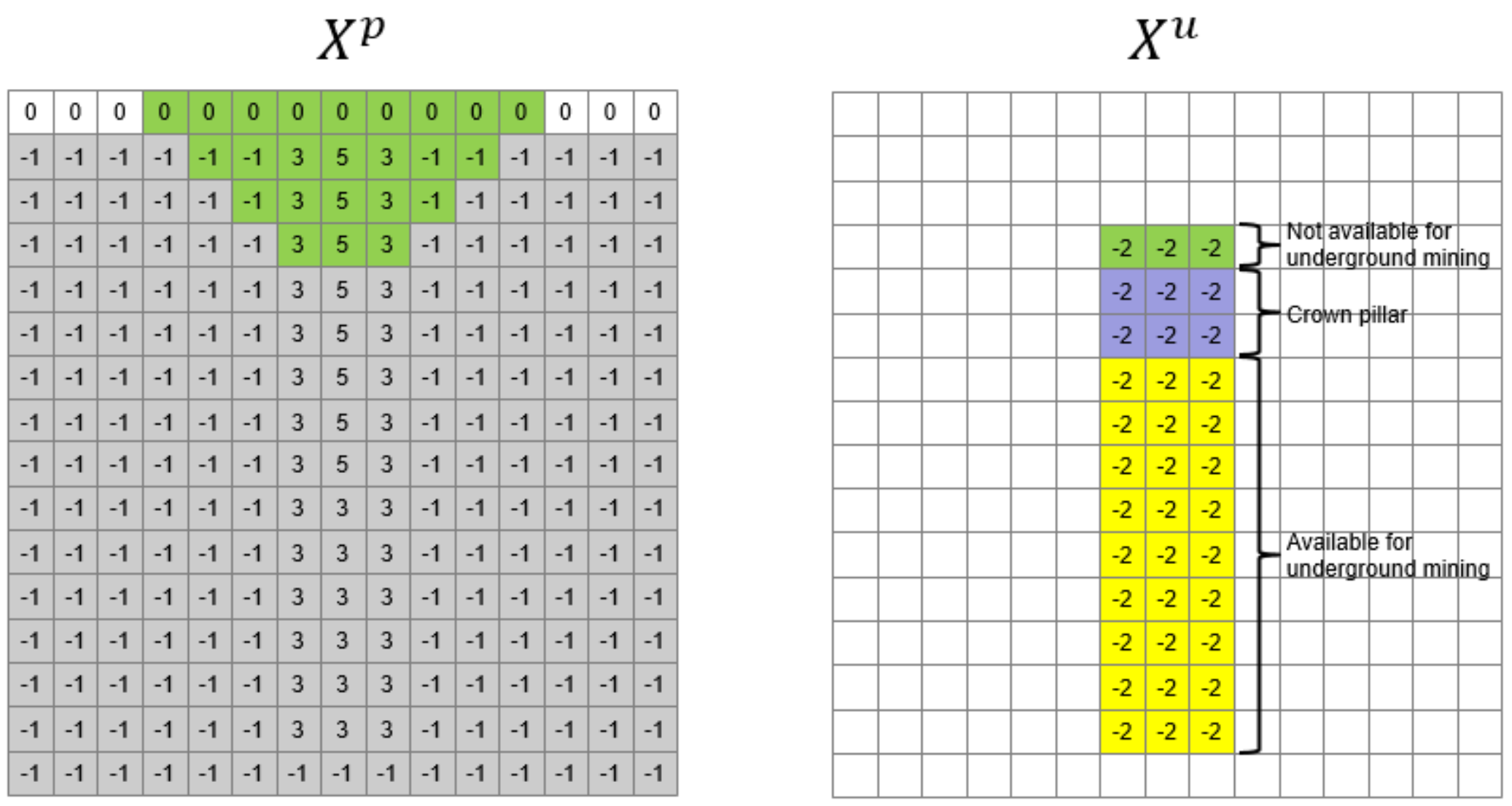

*Figure 9: Arcs $a_i^D \in D$ (not shown) impose a flat top to the crown (see Figure 8 Left). In this case, the bottom of the pit has been flattened as well.*

## The complete modified model

| | |
|---|---|
| $G = (X, A)$ | A digraph representing the open pit optimisation problem with underground option and allowance for a well-formed crown pillar. |
| $x_i^p \in X^p$ | Vertices corresponding to blocks $i$ potentially mineable by open pit method |
| $x_i^u \in X^u$ | Vertices corresponding to blocks $i$ potentially mineable by underground method |
| $X = X^p \cup X^u$ | |
| $B$ | Arcs $a_k^B = (x_i^p, x_j^p)$, $a_k^B \in B$ mapping $X^p$ into $X^p$. Each arc represents a mining dependency to maintain maximum safe pit slopes. |
| $C$ | Arcs $a_k^C = (x_i^p, x_j^u)$, $a_k^C \in C$ mapping $X^p$ into $X^u$<br>Each arc connects vertex $x_i^p$ to vertex $x_j^u$, for every $x_j^u$ that has a weight defined. $x_j^u$ is a block under $x_i^u$, where the vertical separation is equal to the required thickness of the crown pillar. |
| $D$ | Arcs $a_k^D = (x_i^u, x_j^u)$, $a_k^D \in D$ mapping $X^u$ into $X^u$ forming directed cycles at each level comprised of arcs $(x_i^u, x_j^u)$, $x_i^u \neq x_j^u$ . Each cycle represents the required shape for the crown at a given level. |
| $A = B \cup C \cup D$ | |
| $m_i^p = v_i^p - c_i^p \quad \forall x_i^p \in X^p$ | The weight of vertex $x_i^p$ representing the open pit value of the corresponding block. |
| $m_i^u = -(v_i^u - c_i^u) \quad \forall x_i^u \in X^u$ | The weight of vertex $x_i^u$ representing the negative underground value of the corresponding block (the opportunity cost). |
| $G_Y = (Y, A_y)$ | A subgraph of $G$ defined by a set of vertices $Y \subset X$ and containing all the arcs that connect vertices of $Y$ in $G$. |

The problem is to find a closed subgraph $G_Y$ of $G$ that maximises $M_Y$

$$Max\ M_Y = \sum_{i:\, x_i^p \in Y} m_i^p + \sum_{j:\, x_j^u \in Y} m_j^u$$

The results are interpreted as:

| | |
|---|---|
| $x_i^p \in Y$ | Blocks in the optimal open pit. |
| $x_i^u \in Y$ | Blocks not available for underground mining (including blocks in the crown pillar). |
| $P = \{x_i^u \in X^u \cap Y : x_i^p \notin Y\}$ | Blocks in the crown pillar. |
| $X^u \backslash Y$ | Blocks available for underground mining. |

# DEMONSTRATION OF THE NEW OPPORTUNITY COST APPROACH

We demonstrate the new opportunity cost approach through the unconventional use of the Geovia Whittle software. The data we used is a modified version of the well-known Marvin data set.

The original Marvin model comprises around 60 000 blocks of air, undefined waste and mineralised material. To modify the model we first extended the orebody downwards to simulate deep ore potentially minable by open cut or underground method. We then copied this extended model and modified it to create the second underground model. Elevations of the modified Marvin open cut and underground models are shown in Figure 10. The original Marvin is in the top left square. The underground model is to the right. To simulate the case in which not all material is available for underground mining, we removed near-surface oxide material and imposed further arbitrary spatial constraints. These orebodies do not look natural, but their complexity and size are more than adequate to test the new methods. The orebody in this example is contiguous, but the methods will work as well with fragmented orebodies. The models together contain around 250 000 blocks.

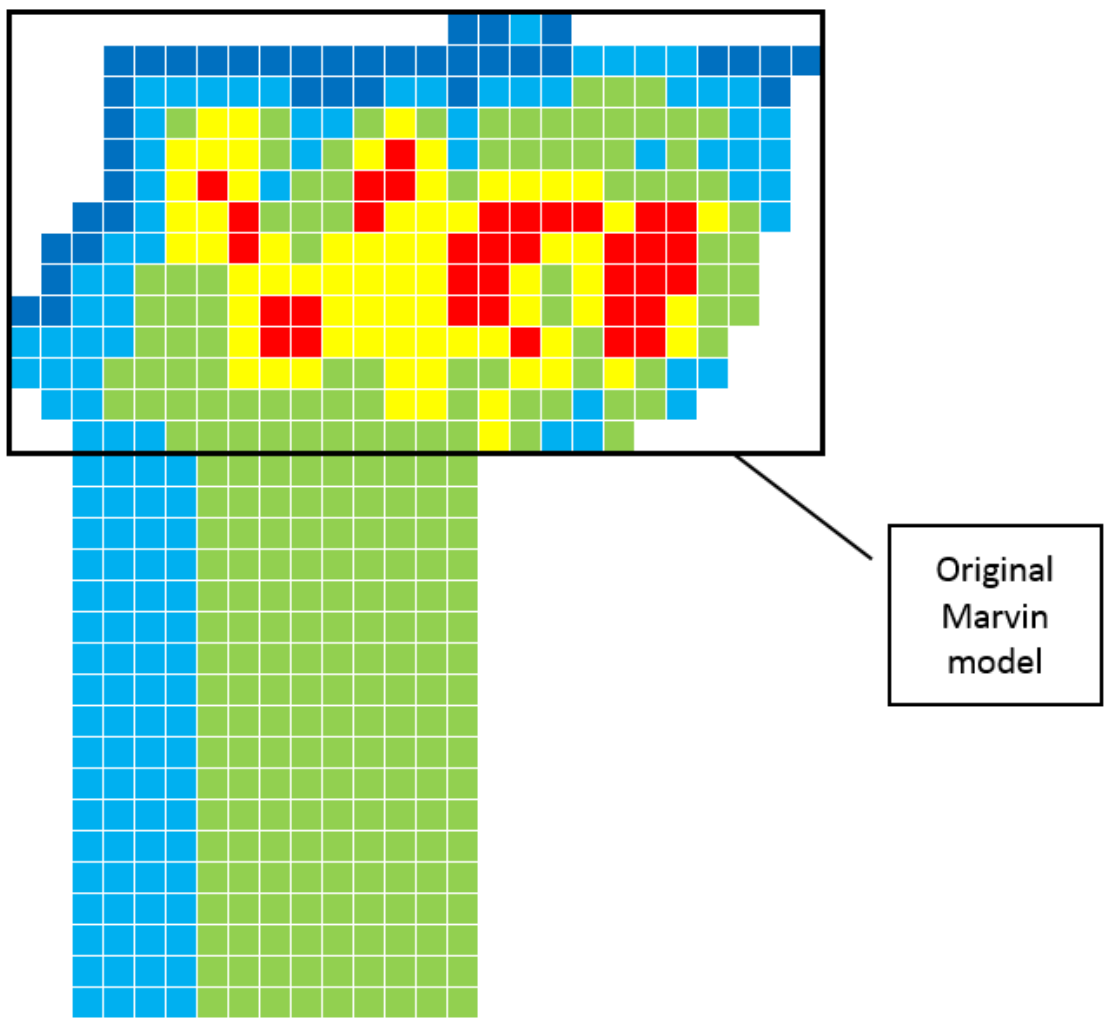

*Figure 10: The original Marvin is in the top left square. The underground model is to the right.*

## Demonstration 1 – Simple Crown Pillar

In this demonstration, we include the first and second modifications to the opportunity cost approach. In order to generate the $B$ arcs, we used the Whittle Slope Set function, with the standard Marvin pit slope settings, but restricted to the open pit model. We created the $C$ arcs in Microsoft Excel and imported them into Whittle.

The result is shown in Figure 11. The black line on the left is the pit. The black outlines on the right are also included in the results, but represent areas that are no longer available for underground mining. The crown pillar is the correct thickness but has an impractical shape. Demonstration 2 shows how a further modification can resolve this issue.

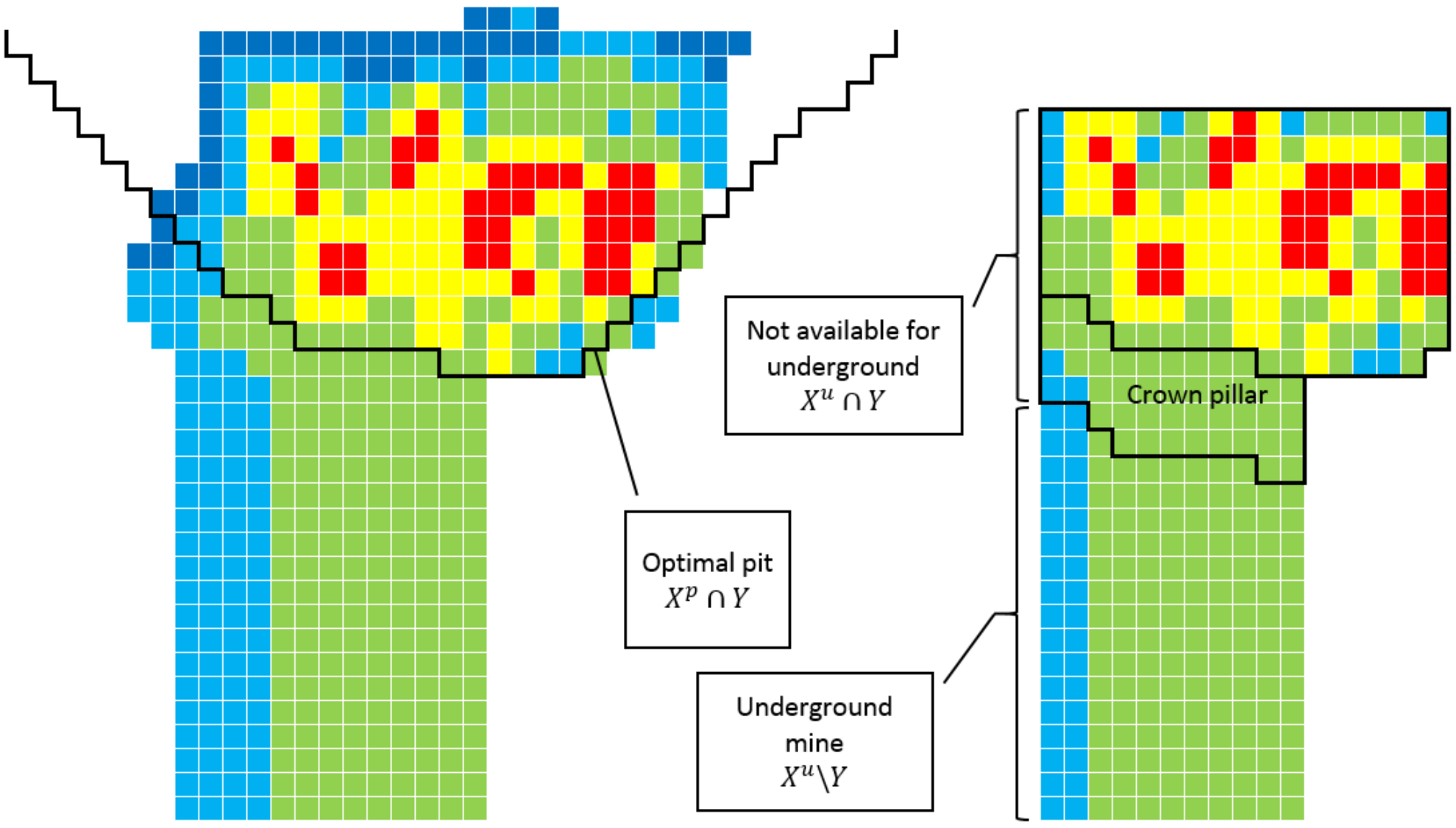

*Figure 11: The black line on the left is the pit. The black outlines on the right are also included in the results, but represent areas that are no longer available for underground mining.*

## Demonstration 2 – Well-formed Crown Pillar

In this demonstration, we include all three modifications to the opportunity cost approach. We created arcs $D$ in Microsoft Excel and imported them into Geovia Whittle. These arcs ensure that a well-formed crown pillars are included in the optimisation. The results are shown in Figure 12.

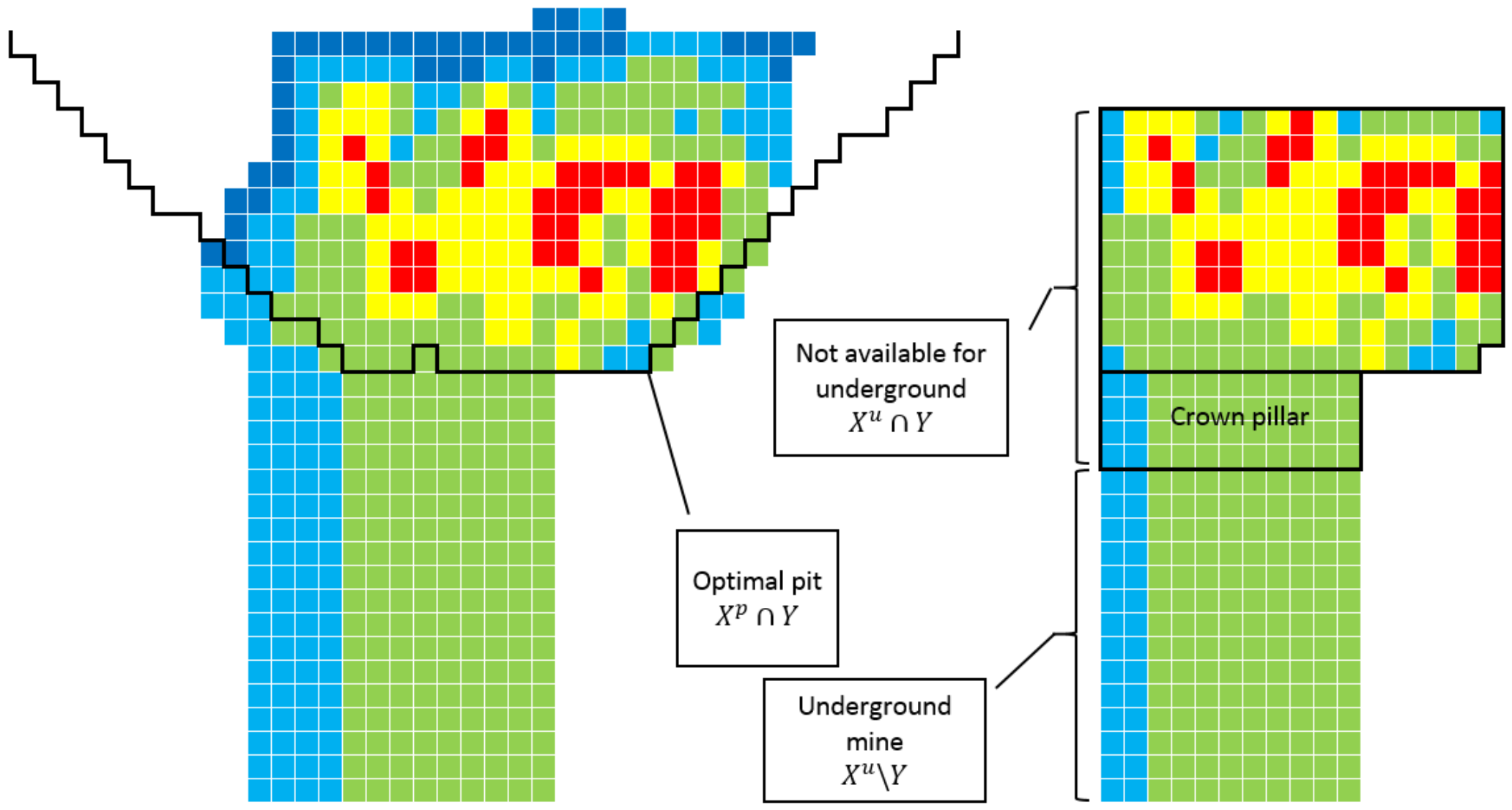

*Figure 12: The crown pillar is the correct thickness and has a shape controlled by the D arcs.*

## Summary Results

Table 1 shows the summary results for the key trials. The first important finding is that the first modification produces identical results to the standard opportunity cost method. With the second modification, we see the impact of a crown pillar and as expected, the value has reduced, since the crown pillar requires that a significant amount of material is left in the ground. Of interest is that the optimal pit is deeper with the crown pillar than without, an outcome particular to this set of data: By pushing the crown pillar down, it sterilises less mineralised material since the orebody narrows. With the third modification we introduced a well formed crown pillar. This additional constraint reduced the value of the mine further as expected.

*Table 1: Summary results from the key trials.*

| **Case**<br>**(Processing time)** | **Arcs** | **Mines** | **Value** |
|---|---|---|---|
| Marvin with vertical extension<br>124 440 blocks<br>(8 sec) | B: 3 296 908<br>C: 0<br>D: 0 | Pit Ore: 350mT<br>Pit Depth: 630m<br>UG Levels: n/a | Open Pit: $1 483m<br>Underground: n/a<br>Total: $1 483m |
| Standard opportunity cost<br>124 440 blocks<br>(8 sec) | B: 3 296 908<br>C: 0<br>D: 0 | Pit Ore: 237mT<br>Pit Depth: 510m<br>UG Levels: 16 | Open Pit: $1 230m<br>Underground: $830m<br>Total: $2 060m |
| First mod. - Alternate opp. cost<br>248 880 blocks<br>(11 sec) | B: 6 733 168<br>C: 13 524<br>D: 0 | Pit Ore: 237mT<br>Pit Depth: 510m<br>UG Levels: 16 | Open Pit: $1 230m<br>Underground: $830m<br>Total: $2 060m |
| Second mod. - Simple crown pillar<br>248 880 blocks<br>(15 sec) | B: 6 733 168<br>C: 13 524<br>D: 0 | Pit Ore: 295mT<br>Pit Depth: 540m<br>UG Levels: 11~19 | Open Pit: $1 396m<br>Underground: $483m<br>Total: $1 879m |
| Third mod. - Well-formed crown pillar<br>248 880 blocks<br>(14 sec) | B: 6 733 168<br>C: 13 524<br>D: 12 075 | Pit Ore: 307mT<br>Pit Depth: 450m<br>UG Levels: 14 | Open Pit: $1 402m<br>Underground: $452m<br>Total: $1 858m |

# PRACTICAL CONSIDERATIONS

The LG algorithm was designed to optimise an open pit outline in order to maximise undiscounted cashflows. Mine planners are most commonly interested in maximising NPV rather than undiscounted value but as discussed in the Introduction, techniques have been developed to enable the LG Algorithm to be used in the pursuit of that objective. In this paper we have developed a modelling approach that extends the use of the LG Algorithm to optimise the pit outline, whilst taking into account the opportunity value of an underground mine and the requirement for a well-formed crown pillar above it. This new method maximises undiscounted cashflows, but as for the regular application of the LG Algorithm, we believe that it can be used to contribute to the maximisation of NPV.

As part of the development of this paper, we solicited a discussion on social media (Whittle et al. 2015). The following is a summary of observations made by contributors:

- Based on a number of actual cases, value (measured as NPV) can be relatively insensitive to a +/- 50 metres change from the optimal transition point.
- Considering the above observation, it is useful to test value sensitivity in each case. This may lead to the identification of options (e.g. several candidate transition depths) all with similar NPVs. The final selection may then proceed on the basis of criteria not directly accounted for in the NPV calculation. For example earlier or later transition may provide for better outcomes in terms of safety, capital scheduling, skills utilisation or production continuity.
- Decisions are always made in conditions of uncertainty, including (but not limited to) price and geological uncertainty. It is useful therefore to use a method to incorporate that uncertainty into the analysis. The development of multiple scenarios is one approach.
- It is important to ensure that cut-off optimisation is undertaken correctly. For example, the failure to take into account the future underground mining, could lead to inappropriate processing of low grade material in the open pit and a loss of value.

We believe in addition these factors should be considered in any application of the approaches developed in this paper:

- Underground mining sequence – In many cases, underground mining sequences commence at the top of the orebody. That material will be mined soon after the completion of the open pit mine. To the extent that differences in timing affects discounting, these effects ought to be easily dealt with by discounting underground values relative to expected delay from the completion of open pit mining. Other underground mining methods (for example block caving) draw material from deep draw points and ore in the transition zone will not be extracted and processed for many years or even decades after the completion of the open pit. In other words, the opportunity cost ought to be heavily discounted and accordingly, will have little or no impact on the open pit outline.
- Iteration – The LG Algorithm-based opportunity cost approach uses a critical assumption: that material in the transition zone will be mined by underground method, if not mined by open cut method. What that practically means is that planners needs to start with an underground mine plan (perhaps a conceptual one), run the open pit optimisation and then consider what is left for the underground. It may be that economic assumptions, especially the calculation of the underground mining extraction cost $c_i^u$ (See the Conventional Opportunity Cost Approach section) are no longer valid for the smaller mining envelope. If this is the case, then a new underground plan should be developed and the open pit optimisation should be run again. We believe that this will be a convergent process, but it is beyond the scope of this paper to discuss this issue in detail.

# FURTHER RESEARCH

As discussed in the Current Approaches section, there are three general approaches to the transition problem, each with their pros and cons. This paper has focused on a new opportunity cost approach. However, there remain gaps in our approach, most particularly the simultaneous optimisation of open cut and underground mines and the direct pursuit of NPV maximisation. Other researchers, often using integer programming approaches, have made some progress in dealing with these gaps but at the expense of model size and accuracy in defining pit slopes, crown pillars and underground stopes.

Our expectation in the long term is that research and development in the field of solving large precedence constrained production scheduling problems will eventually advance to the level that will allow many of the gaps to be bridged.

## CONCLUSIONS

We present in this paper a review of current approaches to the open cut to underground transition problem and develop improvements to the existing opportunity cost approach. The improvements allow the proper consideration of a well-formed crown pillar when optimising the open pit outline above a potential underground mine. The new opportunity cost approach can be implemented with the unconventional application of widely available optimisation tools, as demonstrated with the Geovia Whittle software.

## ACKNOWLEDGEMENT

The authors gratefully acknowledge the support of Dassault Systemes Geovia, especially their provision of an academic license for the Whittle Strategic Mine Planning software. Without this support, it would not have been possible to test and implement the methods developed in this paper.